\documentclass{amsart}
\usepackage{amssymb, amsmath, amsfonts, amscd}

\input xypic

\theoremstyle{plain}
\newtheorem{theorem}{Theorem}[section]
\newtheorem{corollary}[theorem]{Corollary}
\newtheorem{lemma}[theorem]{Lemma}
\newtheorem{proposition}[theorem]{Proposition}

\newtheorem{example}[theorem]{Example}

\theoremstyle{definition}
\newtheorem{definition}[theorem]{Definition}

\theoremstyle{remark}

\numberwithin{equation}{theorem}

\newcommand{\J}{\mathcal{J}}
\newcommand{\A}{\operatorname{A}}
\newcommand{\K}{\operatorname{K}}

\newcommand{\E}{\mathcal{E}}
\renewcommand{\O}{\mathcal{O} }

\renewcommand{\P}{\mathbf{P} }

\newcommand{\Hom}{\operatorname{Hom} }

\newcommand{\End}{\operatorname{End} }
\newcommand{\Ext}{\operatorname{Ext} }

\newcommand{\Diff}{\operatorname{Diff} }
\renewcommand{\H}{\operatorname{H} }
\newcommand{\Proj}{\operatorname{Proj} }

\newcommand{\pone}{\partial_{x} }
\newcommand{\ptwo}{\partial_{y} }
\newcommand{\pthree}{\partial_{z} }
\newcommand{\p}{\partial_{1} }
\newcommand{\pp}{\partial_{2} }
\newcommand{\ppp}{\partial_{3} }

\newcommand{\Spec}{\operatorname{Spec} }

\newcommand{\C}{\mathbf{C} }

\newcommand{\Der}{\operatorname{Der} }

\newcommand{\del}{\partial }
\newcommand{\DR}{\operatorname{DR} }

\newcommand{\Z}{\mathbf{Z}}

\begin{document}

\title{Explcit formulas for algebraic connections on ellipsoid surfaces}


\author{Helge Maakestad }
\address{H\oe gskolen i Bergen} 
\email{Helge.Maakestad@hib.no}



\date{September 2012} 

\begin{abstract} The aim of this paper is to give explicit formulas for algebraic connections on a class of finitely generated 
projective modules $E$ on surface ellipsoids. The connections we construct are non flat with trace of curvature equal to zero hence 
the corresponding Chern class in Lie-Rinehart cohomology is zero. We give a new proof of a formula of the curvature of a connection expressing it 
in terms of the \emph{fundamental matrix} of the module $E$.  In another paper on the subject we introduced a new characteristic classs $c(E) \in
 \operatorname{Ext}^1(L, \End_A(E))$  with the property that
$c(E)$ is trivial iff $E$ has a flat connection. Hence the methods and results in this paper proves that the class $c(E)$ is non trivial in general.
\end{abstract}

\maketitle
\tableofcontents

\section{Introduction}

In general if $A$ is a commutative unital ring and $E$ a finite rank projective $A$-module that is stably trivial, it follows there is a free $A$-module $F=A^m$ and an isomorphism $E\oplus F \cong A^d$. Hence
if $\K(A)$ is the Grothendieck group of finite rank projective $A$-modules, it follows there is an equality $[E]=[A]^{d-m}$ in $\K(A)$. From this it follows all Chern classes $c_i(E) \in \K(A)$ of $E$ in the 
Grothendieck group are trivial. This is true for any theory of Chern classes satisfying the Whitney sum formula. In this paper we show that "most" finite rank projective $A$-modules $E$ have no flat connection. For any $A$-module $E$ there is the Atiyah class

\[ a(E) \in \operatorname{Ext}^1(E, E\otimes \Omega^1_A),  \]

which is zero iff $E$ has a connection $\nabla$. Since $E$ is projective it follows the Atiyah class $a(E)=0$ is zero, hence $E$ always have a connection $\nabla: E \rightarrow E\otimes \Omega^1_A$. The connection $\nabla$ is not unique in general. In this paper we prove that most choices of such a connection $\nabla$ are non-flat by giving an explicit formula for the curvature $R_{\nabla}$ of $\nabla$ in terms of an idempotent endomorphism $\phi$ defining $E$. This formula proves that "most" projective modules $E$ have no flat connection. In another paper on the subject (see \cite{maa145}) I introduce for any Lie-Rinehart algebra $L$ and any $L$-connection $(E, \nabla)$, a pointed non-abelian cohomology set $\operatorname{Ext}^1(L, \End_A(E))$. The connection $(E,\nabla)$ gives rise to a cohomology class 

\[ c(E) \in \operatorname{Ext}^1(L, \End_A(E)),  \]

with the property that $c(E)$ is trivial if and only if $E$ has a flat connection. 
The class $c(E)$ is independent of choice of connection and it gives a class that is non-trivial in general. In particular we get a new characteristic class that is non-trivial for stably trivial projective modules. The class $c(E)$ is defined for an arbitrary module (in particular an arbitrary projective module) on an arbitrary commutative unital ring, and it refines all Chern classes. Hence the methods and results  introduced in this paper proves that the class $c(E)$ introduced in \cite{maa145}  is non-trivial in general.


Associated to a finitely generated projective module $E$ over a commutative unital ring $A$ there is the notion of a projective basis for $E$.
This is a set of generators $e_1,\ldots, e_n$ of $E$ and $x_1,\ldots ,x_n$ of $E^*$ satisfying a compatibility relation. A projective 
basis gives rise to a connection 
\[ \nabla: \Der(A)\rightarrow \Diff^1(E) \]
on the module $E$. The curvature of the connection $\nabla$ satisfies the following formula (see Corollary \ref{lieproduct}): Let $\delta, \eta$ 
be two derivations of $A$. Then there is an equality

\begin{align}
&\label{curvature} R_\nabla(\delta \wedge \eta)=[ \delta(M), \eta(M)] 
\end{align} 
 
in $\End_A(E)$ where $M$ is the fundamental matrix associated to the projective basis $e_i,x_j$. The Chern character and Chern 
classes of $E$ may be defined using
the curvature of a connection on $E$ hence formula \ref{curvature} may be used to describe the Chern character and Chern classes 
of $E$ with values in  algebraic deRham cohomology of $A$. 

If $A$ is a finitely generated algebra over the complex numbers 
it is well known (see \cite{deligne}, \cite{grothendieck} and \cite{hartshorne})  that the algebraic deRham cohomology $\H^*_{\DR}(A)$ 
calculates singular cohomology $\H^*(X_\C, \C)$ 
of the underlying complex variety $X_\C$ of $\Spec(A)$. It follows formula \ref{curvature} might give information on 
the topology of $X_\C$. 

The Chern character
\begin{align}
&\label{chern} Ch: \K_0(A)\rightarrow \H^*_{\DR}(A)
\end{align}
associates to each isomorphism class $[E]$ of a finitely generated projective $A$-module $E$ an element $Ch([E])$ in $H^*_{\DR}(A)$. 
The associated graded group of $\K_0(A)$ with respect to the gamma filtration is isomorphic to the Chow group $\A_*(\Spec(A))$ of 
$\Spec(A)$. Hence the Chern character \ref{chern} gives information on the cycle map
\[ \gamma: \A_*(\Spec(A)) \rightarrow \H^*(X_\C, \C) .\]
It could be formula \ref{curvature} could be used to study algebraic cycles on $\Spec(A)$ and the cycle map $\gamma$.

The calculations of the fundamental matrix $M$, its associated connection $\nabla$ and the curvature $R_\nabla$ may in the case of a 
finitely generated algebra $A$ over a field be done using Groebner bases. 

As a result we calculate in Theorem \ref{main} explicit formulas for non flat algebraic connections
\[ \nabla: \Der_\Z(A) \rightarrow \End_\Z(\Omega) \]
where $\Omega$ is the module of differentials on the ellipsoid surface $x^p+y^q+z^r-1$. The Chern character $Ch(\Omega)$ is zero
(see Lemma \ref{zero} )

We also give explicit examples of non flat algebraic connections 
\[ \nabla: \Der_{\mathbb{C}}(A)\rightarrow \End_A(L) \]
where $L$ is a projective rank one module on the complex two sphere.
We prove the trace of the curvature $R_\nabla$ is non zero (see Proposition \ref{nonzerotrace}).
The associated Chern class $c_1(L)$ in $\H^2_{\DR}(A)$ is non zero 
hence the module $L$ has no flat algebraic connections (see Corollary  \ref{noflat} ).

We also discuss a possible relationship to an old conjecture on existence of flat connections 
on complex vector bundles on complex projective manifolds (see Example \ref{existence}).

\section{Differential operators and idempotents}

In this section we consider $\phi$-connections on free modules and finitely generated projective modules and use this
notion to prove a formula for the curvature $\nabla$ of a connection defined in terms of a projective basis. This formula is 
given in terms of the Lie-product of two matrices.

Let in the following $A$ be a fixed commutative $\Z$-algebra where $\Z$ is a commutative unital ring
 and let $F=A\{u_1,\ldots, u_n\}$ be a free $A$-module of rank $n$.
Let $E$ be a finitely generated projective module and let 
\[ 0\rightarrow K \rightarrow F \rightarrow^p E\rightarrow E \rightarrow 0 \]
be an exact sequence. Since $E$ is projective there is an $A$-linear  section $s$ of $p$. Let in the following $\phi=s\circ p$ and
$\psi=I_n-\phi$ where $I_n$ is the $n\times n$ identity matrix. It follows $\phi, \psi\in \End_A(F)$.

\begin{lemma} The following holds: $\phi^2=\phi$, $\psi^2=\psi$ and $\phi + \psi=I$. Moreover $ker(\phi)=K$ and $ker(\psi)\cong E$.
\end{lemma}
\begin{proof} It is clear thet $\phi^2=\phi$. We moreover get
\[\psi^2= (I_n-\phi)(I_n-\phi)=I_n^2-I_n\circ \phi -\phi\circ I_n -\phi^2=\]
\[ I_n -2\phi +\phi =I_n-\phi=\psi.\]
Clearly $\phi+\psi=I_n$ and the claim follows.
We claim the following equality: $ker(\phi)=K$.
Assume $x\in ker(\phi)$. It follows $s(p(x))=0$ hence $p(s(p(x))=p(x)=0$ hence $x\in ker(p)=K$. It follows 
$ker(\phi)\subseteq K$. Assume $x\in K=ker(p)$. It follows $s(p(x))=s(0)=0$ hence $x\in ker(\phi)$ and the claim follows.
Consider the induced map $p:ker(\psi)\rightarrow E$. We claim this map is an isomorphism. 
Injectivity: Assume $x\in ker(\psi)$. It follows $\phi(x)=x$. Assume $p(x)=0$. It follows $x\in ker(p)=ker(\phi)$. It follows
$x=\phi(x)=0$ hence $p|_{\psi}$ is injective. We prove $p|_{\psi}$ is surjective. Assume $e\in E$. It follows $p(s(e))=e$. Let
$x=s(e)$. We get
\[ \phi(x)=\phi(s(e))=s(p(s(e))=s(e)=x.\]
It follows $\phi(x)=x$ hence $x\in ker(\psi)$. It follows $p|_{ker(\psi)}$ is surjective and the claim follows.
\end{proof}

\begin{lemma} The canonical map $f:ker(\phi)\oplus ker(\psi)\rightarrow F$ defined by
\[ f(x,y)=x+y \]
is an isomorphism of $A$-modules.
\end{lemma}
\begin{proof} Assume $f(x,y)=x+y=0$. It follows $x=-y$. It follows $-x=y\in ker(\phi)$. We get
\[ 0=\phi(-x)=-\phi(x)=-x \]
hence $x=0$. It follows $(x,y)=(-y,y)=(0,0)=0$ and the map $f$ is injective.
Surjectivity: Assume $x\in F$. It follows $x-s(p(x))\in K$ and $p(x)\in E$.  It follows $s(p(x))\in ker(\psi)$:
\[ \psi(s(p(x)))=(I-\phi)(s(p(x)))=s(p(x))-s(p(s(p(x))))=s(p(x))-s(p(x))=0.\]
It follows $(s(p(x)), x-s(p(x))\in ker(\psi)\oplus ker(\phi)$. Moreover
\[ f(s(p(x)), x-s(p(x)) )=x.\]
The claim follows.
\end{proof}

\begin{definition} Let $\phi\in \End_A(E)$ be an endomorphism. A map of $A$-modules
\[ \nabla: \Der_{\Z}(A)\rightarrow \End_{\Z}(E) \]
satisfying
\[ \nabla(\delta)(ae)=a\nabla(\delta)(e)+\delta(a)\phi(e) \]
is called a $\phi$-connection on $E$. 
\end{definition}

Let $E,F$ be left $A$-modules and define $\Diff^0(E,F)=\Hom_A(E,F)$ Let $k\geq 1$ be an integer and define
\[ \Diff^k(E,F)=\{ \partial \in \Hom_\Z(E,F): [\partial, \phi_A]\in \Diff^{k-1}(E,F)\text{ for all }a\in A\}.\]
Let $\Diff(E,F)=\cup_{k\geq 0}\Diff^k(E,F)$. Let $\Diff(E)=\Diff(E,E)$. It follows $\Diff(E)$ is an associative sub ring of $\End_\Z(E)$.
called the \emph{ring of differential operators on $E$}. The module $E$ is a $\Diff(A)$-module if and only if there is a ring homomorphism
\[ \rho: \Diff(A) \rightarrow \Diff(E).\]
If $A$ is a polynomial ring over a field of characteristic zero it follows 
$\Diff(A)$ is the Weyl algebra.

\begin{lemma} A $\phi$-connection induce a map
\[ \nabla: \Der_{\Z}(A)\rightarrow \Diff^1(E).\]
\end{lemma}
\begin{proof} We need to check that for every derivation $\delta \in \Der_{\Z}(A)$ it follows
$\nabla(\delta)\in \Diff^1(E)$. Let $a\in A$ and consider the following endomorphism of $E$:
\[ [\nabla(\delta),a](e)=\nabla(\delta)(ae)-a\nabla(\delta)(e)=\]
\[ a\nabla(\delta)(e)+\delta(a)\phi(e)-a\nabla(\delta)(e)=\delta(a)\phi(e).\]
It follows $[\nabla(\delta),a]=\delta(a)\phi$. Hence $[\nabla(\delta),a]\in \End_A(E)$  for every
$a\in A$. The claim follows.
\end{proof}

Note: If $\phi$ is the zero endomorphism a $\phi$-connection is a morphism
\[ \nabla: \Der_{\Z}(A)\rightarrow \End_A(E).\]
A $0$-connection is also called a \emph{potential}.
If $\phi$ is the identity morphism it follows 
\[\nabla:\Der_{\Z}(A)\rightarrow \End_{\Z}(E) \]
is an ordinary connection on $E$. 

Note that the sum of a $\phi$-connection and a $\psi$-connection is a $\phi+\psi$-connetion. It follows that the set of
connections on $E$ is a torsor on the set of potentials: Given a connection $\nabla$ and a potential $\phi$ it follows 
$\nabla +\phi$ is a connection. All connections on $E$ arise in this way.

The free $A$-module $F$ has a canonical flat connection defined as follows:
Let $\delta$ be a dervation of $A$ and defined
\[ D_{\delta}:F\rightarrow F \]
by
\[ D_{\delta}(\sum a_i u_i)=\sum \delta(a_i)u_i.\]
It follows we get a map
\[ D: \Der_{\Z}(A)\rightarrow \End_{\Z}(F) .\]
Define for every derivation $\delta$ the following endomorphism of $F$:
\[ \nabla(\delta)(x)=D_{\delta}(\phi(x)) \]
where $\phi \in \End_A(F)$ is the idempotent defined above.

\begin{lemma} The map $\nabla$ is a $\phi$-connection. 
\end{lemma}
\begin{proof} We get the following calculation:
\[ \nabla(\delta)(ae)=D_{\delta}(\phi(ae))=D_{\delta}(a\phi(e))=\]
\[ aD_{\delta}(\phi(e))+\delta(a)\phi(e)=a\nabla(\delta)(e)+\delta(a)\phi(e).\]
The claim follows.
\end{proof}

\begin{lemma} The $\phi$-connection $\nabla$ induce a connection $\del$ on $E$ via $\del(\delta)(e)=p(\nabla(\delta)(e))$.
\end{lemma}
\begin{proof} Assume $p(x)=p(y)=e\in E$. Define $\\del(\delta)(e)=p(\nabla(\delta)(y))$. We show this is a well defined connection 
on $E$. We get $x-y=k\in K=ker(p)$.
We get
\[ \del(\delta)(e)=p(\nabla(\delta)(y+h)=p(\nabla(\delta)(y))+p(\nabla(\delta)(h))=\]
\[ p(\nabla(\delta)(y))+p(D_{\delta}(s(p(h))))=p(\nabla(\delta)(y).\]
It follows $\del$ is well defined.
Let $a\in A$ and $p(ay)=ap(y)=ae$. We get
\[ \del(\delta)(ae)=p(\nabla(\delta)(ay)=\]
\[ p(a\nabla(\delta)(y)+\delta(a)\phi(y))=ap(\nabla(\delta)(y))+\delta(a)p(\phi(y))=\]
\[ a\del(\delta)(e)+\delta(a)p(s(p(y))=a\del(\delta)(e)+\delta(a)e.\]
The claim follows.
\end{proof}

Let $e_i=p(u_i)$ be generators of the $A$-module $E$. Let $u_i^*:F\rightarrow A$ be dual to the free basis $u_1,\ldots , u_n$. 
Let $x_i=u_i^*\circ s$. It follows the set $e_i, x_j$ form a projective basis for $E$ in the sense of \cite{maa14}. 
Let

\[ 
M=\begin{pmatrix} x_1(e_1) & \cdots   & x_1(e_n) \\
                   \vdots  &  \cdots  & \vdots   \\
                  x_n(e_1) & \cdots   & x_n(e_n)
\end{pmatrix}
\]
be the fundamental matrix of $E$ with respect to the pair $p,s$.
Since $M$ is an $n\times n$-matrix with coefficients in $A$ it follows we get an $A$-linear map $M:F\rightarrow F$.

\begin{lemma} There is an equality of maps $M=\phi$.
\end{lemma}
\begin{proof} Let $x=a_1u_1+\cdots +a_nu_n\in F$ be an element.  It follows $p(x)=\sum a_ie_i$.
We get the following calculation: $\phi(x)= s(p(x))=s(\sum a_i e_i)=\sum a_i s(e_i)$. Hence
The matrix of $\phi$ in the basis $u_j$ has $ s(e_j)$ as column vectors. We see that 
\[ s(e_j)= b_1u_1+\cdots +b_nu_n=\]
\[ u_1^*(s(e_j))u_1+\cdots + u_n^*(s(e_j))u_n=\]
\[x_1(e_j)u_1+\cdots +x_n(e_j)u_n.\] 
It follows $s(e_j)$ is the vector
\[
s(e_j)=
\begin{pmatrix} x_1(e_j) \\
              \vdots \\
              x_n(e_j)
\end{pmatrix}.
\]
It follows the matrix of $\phi$ in the basis $u_j$ is the matrix $M$ and the claim follows.
\end{proof}

\begin{lemma} The following formula holds:
\[
\phi(x)=
\begin{pmatrix}x_1(p(x)) \\
          \vdots \\
               x_n(p(x))
\end{pmatrix}.
\]
\end{lemma}
\begin{proof} The proof is left to the reader.
\end{proof}

Consider the following endomorphism of $F$:
\[ g=[D_{\delta}, \phi]:F\rightarrow F.\]
Let $\delta(M)$ be the $n\times n$-matrix with $\delta(x_i(e_j))$ as entries.

\begin{lemma} The endomorphism $[D_{\delta}, \phi]$ is $A$-linear and its matrix in the basis $u_j$ is $\delta(M)$.
\end{lemma}
\begin{proof} The proof is left to the reader.
\end{proof}

By definition $\nabla(\delta)=D_{\delta}\circ \phi$. It follows $\nabla:\Der_{\Z}(A)\rightarrow \End_{\Z}(E)$ is a
$\phi$-connection. Let $\del(\delta)=p\circ \nabla(\delta)\circ s:E\rightarrow E$. It follows 
$\del(\delta)=p\circ D_{\delta}\circ s$ and from the calculations above it follows $\del$ is a connection on $E$.
For a pair of derivations $\delta, \eta$ the curvature $R_\nabla(\delta \wedge \eta):E\rightarrow E$ is an $A$-linear
endomorphism. 

\begin{lemma} There is an  equality $R_\del(\delta\wedge \eta)=p\circ R_\nabla(\delta \wedge \eta) \circ s$ of maps.
\end{lemma}
\begin{proof} By definition it follows $\del(\delta)=p\circ \nabla(\delta)\circ s$. We get
\[ R_\del(\delta \wedge \eta)=[\del(\delta),\del(\eta)]-\del([\delta,\eta]).\]
Since $p\circ s =id$ this equals
\[ p\circ \nabla(\delta)\circ \nabla(\eta)\circ s - p\circ \nabla(\eta)\circ \nabla(\delta)\circ p - p\circ \nabla([\delta,\eta])\circ s=\]
\[ p\circ(\nabla(\delta)\circ \nabla(\eta)-\nabla(\eta)\circ \nabla(\delta)-\nabla([\delta,\eta]))\circ s=p\circ R_\nabla(\delta \wedge \eta)\circ s\]
and the claim follows.
\end{proof}

\begin{lemma} The following formula holds:
\[ \del(\delta)(e)=\delta(x_1(e))e_1+\cdots +\delta(x_n(e))e_n.\]
\end{lemma}
\begin{proof} By definition it follows $\del(\delta)(e)=p(\nabla(\delta)(s(e)))=p(D_{\delta}(\phi(s(e))))$.
By the above results it follows $\phi(s(e))$ equals the vector
\[
\begin{pmatrix} x_1(p(s(e)) \\
                 \vdots  \\
               x_n(p(s(e))
\end{pmatrix}.
\]
Since $p(s(e))=e$ we get
\[
\phi(s(e))=\begin{pmatrix} x_1(e) \\
                 \vdots  \\
               x_n(e)
\end{pmatrix}.
\]
We get $\nabla(\delta)(e)=p(D_\delta(\phi(s(e)))=\delta(x_1(e))e_1+\cdots +\delta(x_n(e))e_n$ and the claim follows.
\end{proof}

Let for any derivation $\delta$ $A_\delta$ be the endomorphism
\[ A_\delta:D_\delta+\delta(M):F\rightarrow F.\]

\begin{lemma} The following formula holds:
\[ [A_\delta, A_\eta]-A_{[\delta, \eta]}=[\delta(M),\eta(M)] .\]
\end{lemma}
\begin{proof} This is a straight forward calculation: 
\[ [A_\delta, A_\eta]-A_{[\delta, \eta]}=\]
\[ (D_\delta+\delta(M))(D_\eta+\eta(M))-(D_\eta+\eta(M))(D_\delta+\delta(M))-D_{[\delta,\eta]}-[\delta,\eta](M)=\]
\[D_\delta D_\eta+ D_\delta \eta(M)+\delta(M)D_\eta +\delta((M)\eta(M) \]
\[ -D_\eta D_\delta -D_\eta \delta(M)-\eta(M)D_\delta-\eta(M)\delta((M)-D_{[\delta, \eta]}-[\delta,\eta](M)=\]
\[ [D_\delta, \eta(M)] +[\delta(M), D_\eta]+[\delta(M),\eta(M)]-[\delta,\eta](M)\]
\[ (\delta\circ \eta)(M)-(\eta\circ \delta)(M)+[\delta(M),\eta(M)]-[\delta,\eta](M)=[\delta(M),\eta(M)]\]
and the claim follows.
\end{proof}

\begin{lemma} \label{derivation} There is an equality of maps $p \circ A_\delta=\del(\delta)\circ p$
where $\del(\delta)=p\circ \nabla(\delta)\circ s$ and $\nabla(\delta)=D_\delta\circ \phi$.
\end{lemma}
\begin{proof} We get
\[ p\circ A_\delta=p\circ(D_\delta + \delta(M))=p\circ(D_\delta+ D_\delta\circ \phi -\phi\circ D_\delta)=\]
\[p\circ D_\delta +p\circ D_\delta \circ \phi -p\circ \phi\circ D_\delta=\]
\[p\circ D_\delta +p\circ D_\delta \circ \phi -p\circ D_\delta= p\circ D_\delta \circ \phi=\]
\[ p\circ D_\delta \circ s \circ p=\nabla(\delta)\circ p\]
and the claim follows.
\end{proof}

\begin{theorem} \label{mainlie} There is for every pair of derivations $\delta, \eta$ an equality
\[R_\nabla(\delta\wedge \eta)\circ p=p\circ([A_\delta,A_\eta]-A_{[\delta,\eta]}).\]
\end{theorem}
\begin{proof} Since $p\circ A_\delta=\nabla(\delta)\circ p$ the Theorem follows.
\end{proof}

\begin{corollary} \label{lieproduct} Assume $p(x)=e\in E$. The following formula holds:
\[ R_\nabla(\delta\wedge \eta)(e)=p([\delta(M), \eta(M)](x)).\]
\end{corollary}
\begin{proof} By the above calculation the following holds:
\[ R_\nabla(\delta \wedge \eta)(e)=R_\nabla(\delta \wedge \eta)(p(x))=\]
\[p([A_\delta, A_\eta]-A_{[\delta,\eta]})(x)=p([\delta(M),\eta(M)](x)) \]
and the Corollary follows.
\end{proof}

Let $g:E\rightarrow E$ be an $A$-linear endomorphism and let $\tilde{g}=s\circ g\circ p$ be an endomorphism of $F$.

\begin{lemma} There is an equality $tr(\tilde{g})=tr(g)$ in $A$.
\end{lemma}
\begin{proof} Using the canonical isomorphism
\[ \rho: E^*\otimes_A E\rightarrow \End_A(E) \]
we define the trace
\[ tr:E^*\otimes_A E\rightarrow A\]
by
\[ tr( g\otimes e)=g(e).\]
Let $F=A\{u_1,\ldots ,u_n\}$ and $p:F\rightarrow E$ with $p(u_i)=e_i$. Let $F^*=A\{y_1,\ldots ,y_n\}$ with
$y_i=u_i^*$. Let $s$ be a section of $p$. By definition $x_i=u_i^*\circ s$. Using the isormophism $\rho$ it follows we may write
\[ \tilde{g}=x_1\otimes g(e_1)+\cdots +x_n\otimes g(e_n) \]
in $E^*\otimes_A E$. It follows $\rho(\tilde{g})=g$. We get the following calculation:
\[ tr(g)=tr(x_1\otimes g(e_1)+\cdots + x_n\otimes g(e_n))=\]
\[x_1(g(e_1))+\cdots +x_n(g(e_n)).\]
We moreover get
\[ tr(\tilde{g})=tr(y_1\otimes \tilde{g}(u_1)+\cdots +y_n\otimes \tilde{g}(u_n)=\]
\[ y_1(\tilde{g}(u_1))+\cdots +y_n(\tilde{g}(u_n))=\]
\[u_1^*(s\circ g\circ p(u_1))+\cdots u_n^*(s\circ g\circ p(u_n))=\]
\[ (u_1^*\circ s)(g(e_1)+\cdots +(u_n^* \circ s)(g(e_n))=\]
\[ x_1(g(e_1))+\cdots +x_n(g(e_n))=tr(g)\]
and the claim follows.
\end{proof}

Let 
\[ 0 \rightarrow E' \rightarrow E \rightarrow E" \rightarrow 0 \]
be an exact sequence of finitely generated projective $A$-modules and let $g:E\rightarrow E$ be an $A$-linear endomorphism
with $g(E')\subseteq E'$. It follows $g$ induce an $A$-linear endomorphism $\tilde{g}\in \End_A(E")$.
Let $\Z \subseteq A$ the minimal subring containing the unity of $A$ and let $\Z[t]$ be the polynomial ring in $t$ over $\Z$.
Let $A[t]=A\otimes_{\Z} \Z[t]$ be the ring of polynomials in $t$ with coefficients in $A$. Let $E[t]=E\otimes_A A[t]$. It follows
$E[t]$ in a canonical way is an $A[t]$-module. We get in a cannonical way an exact sequence
\[ 0\rightarrow E'[t]\rightarrow E[t]\rightarrow E"[t]\rightarrow 0\]
of $A[t]$-modules. The endomorphism $g$ induce in a canonical way an endomorphism $g_t:E[t]\rightarrow E[t]$ with 
$g_t(E'[t])\subseteq E'[t]$. Let $I$ be the identity endomorphism of $E$. It follows $G=tI-g_t$ is an endomorphism of $E[t]$ and
since $E$ is locally free of finite rank we may locally take the determinant $det(G)$ of $G$ since $G$ locally is given as a
square matix with coefficients in a commutative ring. This construction glues to give a determinant $P_g(t)=det(G)$ of $G$ which is a polynomial in 
$t$ with coefficients in $A$. Let $e$ be the rank of $E$. It follows $deg(P_g(t))=e$.
We let $P_g(t)$ be the \emph{characteristic polynomial of $g$}. Since $g$ fixes the module $E'$ we get an induced endomorphism $g'$ og $E'$ 
and $g"$ of $E"$. The following holds:
\begin{lemma} \label{char} There is an equality
\[ P_g(t)=P_{g'}(t)P_{g"}(t) \]
of polynomials in $A[t]$.
\end{lemma}
\begin{proof}The proof is left to the reader.
\end{proof}

From Corollary \ref{lieproduct} there is an equality
\[ R_\nabla(\delta\wedge \eta)=p \circ [\delta(M),\eta(M)]\circ s \]
of maps of $A$-modules. The endomorphisms $\delta(M), \eta(M)$ do not induce elements in $\End_A(E)$ hence
the commutator $[\delta(M),\eta(M)]$ is not the commutator of elements in $\End_A(E)$. 
Let 
\[ 0\rightarrow K \rightarrow A^n \rightarrow E \rightarrow 0\]
be the exact sequence with $K\oplus E\cong A^n$.  Let $x\in K\subseteq A^n$. It follows
\[ p([\delta(M),\eta(M)](x)=R_\nabla(\delta\wedge \eta)(p(x))=0 .\]
It follows 
\[ [\delta(M),\eta(M)](x)\in K\]
hence the $A$-linear map $g=[\delta(M),\eta(M)]$ fixes the $A$-module $K$ and induce the map $R_\nabla(\delta\wedge \eta)$ on $E$. 
Let $g_K$ be the induced endomorphism of $g$ on $K$.
It follows from Lemma \ref{char} there is an equality of characterictic polynomials
\[ P_g(t)=P_{g_K}(t)P_{R_\nabla(\delta \wedge \eta)}(t).\]
It follows 
\[ 0=tr(g)=tr(g_K)+tr(R_\nabla(\delta \wedge \eta))\]
hence
\[ tr(R_\nabla(\delta \wedge \eta))=-tr(g_K).\]

\begin{proposition} \label{tracecurv} There is an equality
\[ tr(R_\nabla(\delta \wedge \eta))=-tr(g_K)\]
in $A$
\end{proposition}
\begin{proof} The proof follows from the discussion above.
\end{proof}

\begin{corollary} If $tr(g_K)=0$ it follows $tr(R_\nabla(\delta \wedge \eta))=0$.
\end{corollary}
\begin{proof} This follows from Proposition \ref{tracecurv}.
\end{proof}

The trace of $[\delta(M),\eta(M)]$ on $A^n$ is zero however.

The formula 
\[ R_\nabla(\delta \wedge\eta)=p\circ [\delta(M),\eta(M)]\circ s \]
is an algebraic formula giving  the curvature of $\nabla$ in terms of the matrix $M$. 

\begin{example} The Chern character with values in algebraic DeRham cohomology.\end{example}

The matrix $M$ may in the case where $A$ 
is a finitely generated algebra over a field and $E$ a finitely generated projective $A$-module be calculated using Groebner bases.
One needs to find a set of generators $e_1,\ldots , e_n$ of $E$ and calculate a splitting $s$ of the corresponding surjection
$p:A^n\rightarrow E$ defined by $p(u_i)=e_i$. 

In the case of the Chern character
\[ ch:\K_0(A)\rightarrow \H^*_{\DR}(A) \]
if one knows generators for $\K_0(A)$ as an abelian group one may for each such generator $E$ calculate explicit formulas for a connection
$\nabla$. In such cases one may give a description of the Chern character $ch$ since we have complete control on the curvature $R_\nabla$. 

Let $X=\Spec(A)$.
When $A$ is a finitely generated regular algebra over the complex numbers and $X_{\C}$ is the underlying complex manifold of $X$ it follows
\[ \H^*_{\DR}(A)\cong \H^*(X_\C, \C) \]
where $\H^*(X_\C,\C)$ is singular cohomology of $X_\C$ with $\C$-coefficients. In this case we may describe the image of the Chern character
in singular cohomology and compare this to the image of the cycle map
\[ \gamma: A_*(X)\otimes_\Z \mathbb{Q}\rightarrow \H^*(X_\C, \C) \]
from the Chow group of algebraic cycles. 

\begin{example} An invariant for finitely generated projective modules.\end{example}

Let in the following $E$ be a finitely generated projective $A$-module where $A$ is a commutative unital ring.

\begin{definition}
Let $pr(E)$ be the minimum of the number of elements $n$ in a projective basis $e_1,\ldots , e_n,x_1,\ldots , x_n$ for $E$.
The number $pr(E)$ is the \emph{projective rank} of $E$.
Let $dev(E)=pr(E)-rk(E)$ be the \emph{deviation of} $E$, where $rk(E)$ is the \emph{rank} of $E$ as a projective $A$-module.
\end{definition}

One sees that $dev(E)\geq 0$ for any finitely generated projective module $E$

\begin{proposition} Assume $A$ is a noetherian ring. Then $dev(E)=0$ if and only if $E$ is a free $A$-module.
\end{proposition}
\begin{proof} The proof is left to the reader as an exercise.
\end{proof}

Note: The cotangent bundle $\Omega$ on the real two sphere is a non trivial rank two projective module. One checks that
$dev(\Omega)=1$ in this case.

\section{Algebraic connections on cotangent bundles on ellipsoid surfaces}

The aim of this section is to show that the tecniques introduced in the previous section may be used to give new examples of explicit formulas for 
algebraic connections on a large class of finitely generated projective modules on ellipsoid surfaces. Assume $A$ is a finitely 
generated regular algebra over a field $\Z$ of characteristic zero and $E$ is a finitely generated projective $A$-module. It follows the first order
jet module $\J^1(E)$ is a finitely generated and projective $A$-module. Since $E$ is projective it has a connection and a connection $\nabla$ is 
given by a left splitting $s$ of the canonical map
\[ p:\J^1(E)\rightarrow E \]
of left $A$-modules. It follows a connection $\nabla$ on $E$ is algebraic and may be described in terms of polynomial functions. The aim 
of this section is to give explicit algebraic formulas for a class of connections on modules of differentials on surface ellipsoids.

Given the polynomial $f(x,y,z)=x^p+y^q+z^r-1\in \Z[x,y,z]$ 
let $A=\Z[x,y,z]/f$. Let $X=\Spec(A)$ and let $\Omega=\Omega^1_{X/\Z}$ be the module of Kahler-differentials on $X$. Since $X$ is a smooth
scheme over $\Z$ it follows $\Omega$ is a locally free $A$-module of rank two. In this section we give explicit algebraic fomulas
for algebraic connections on $\Omega$ for any integers $p,q,r\geq 2$ and show these connections are non-flat in general. 

In the paper \cite{maa1} the Kodaira-Spencer map
\[ g: \Der_\Z(A)\rightarrow \Ext^1_A(\Omega,\Omega) \]
and the Kodaira-Spencer class 
\[ a(\Omega)\in \Ext^1_A(\Der_\Z(A), \End_A(\Omega)) \]
was used to give such formulas. As shown in section one of this paper, when $\Omega$ is locally free one can calculate
explicit formulas using the fundamental matrix $M$ of a split surjection $A^n\rightarrow \Omega$. 

All explicit examples of connections
on isolated hypersurface singularities given in the paper \cite{maa1} are flat. The connections calculated in this section are non flat.

Consider the differential $df$ of $f$:
\[ df=px^{p-1}dx+qy^{q-1}dy+rz^{r-1}dz \in A\{ dx,dy,dz\}.\]
By definition $\Omega=A\{dx,dy,dz\}/df$. Let $e_1=\overline{dx}, e_2=\overline{dy}$ and $e_3=\overline{dz}$.
We get a surjection
\[p:A\{u_1,u_2,u_3\}\rightarrow \Omega \]
defined by $p(u_i)=e_i$ for $i=1,2,3$. 
Let 
\[ dF=px^{p-1}u_1+qy^{q-1}u_2+rz^{r-1}u_3 \in A\{u_1,u_2,u_3\}.\]
Make the following definition:
\[ s(e_1)=u_1-\frac{x}{p}dF \]
\[ s(e_2)=u_2-\frac{y}{q}dF \]
and
\[ s(e_3)=u_3-\frac{z}{r}dF.\]

\begin{lemma} The map $s$ is an $A$-linear section of $p$.
\end{lemma}
\begin{proof} The following holds:
\[ 0=s(0)=s(px^{p-1}e_1+qy^{q-1}e_2+rz^{r-1}e_3)=\]
\[ px^{p-1}(   u_1- \frac{x}{p}dF       )+qy^{q-1}(u_2-\frac{y}{q}dF  )+rz^{r-1}( u_3-\frac{z}{r}dF  )=\]
\[ dF-(x^p+y^q+z^r)dF=dF-dF=0\]
It follows $s$ is a well defined $A$-linear section of $p$.
\end{proof}

One calculates
\[ s(e_1)=(1-x^p)u_1-\frac{q}{p}xy^{q-1}u_2 -\frac{r}{p}xz^{r-1}u_3,\]
\[ s(e_2)=  -\frac{p}{q}x^{p-1}yu_1+(1-y^{q})u_2 -\frac{r}{q}yz^{r-1}u_3\]
and
\[  s(e_3)=-\frac{p}{r}x^{p-1}zu_1-\frac{q}{r}y^{q-1}zu_2+(1-z^r)u_3.\]

It follows the matrix $M$ is as follows:
\[
M=
\begin{pmatrix} 1-x^p    &   -\frac{p}{q}x^{p-1}y     &  -\frac{p}{r}x^{p-1}z    \\
              -\frac{q}{p}xy^{q-1}  &   1-y^q    &   -\frac{q}{r}y^{q-1}z        \\
              -\frac{r}{p}xz^{r-1}  &   -\frac{r}{q}yz^{r-1}   &  1-z^r 
\end{pmatrix}.
\]
Let 
\[ \p=f_y\pone -f_x\ptwo =qy^{q-1}\pone -px^{p-1}\ptwo ,\]
\[\pp=f_z\pone -f_x\pthree=rz^{r-1}\pone -px^{p-1}\pthree\]
and
\[ \ppp=f_z\ptwo -f_y\pthree=rz^{r-1}\ptwo -qy^{q-1}\pthree.\]
It follows the derivations $\p,\pp, \ppp$ generate the $A$-module $\Der_\Z(A)$. In the following we give explicit formuls for:
First order differential operators
\[ \nabla(\partial_i)=D_{\partial_i}+\partial_i(M):\Omega \rightarrow \Omega \]
for $i=1,2,3$. 


Let in the following $\partial \in \Der_\Z(A)$ be a derivation and let 
\[ D_\partial:A^3\rightarrow A^3 \]
be the following operator:
\[ D_\partial(\sum a_i u_i)=\sum \partial(a_i)u_i.\]
It follows $D_\partial$ is a first order differential operator on $A^3$. 
Let $\partial(M)$ be the $3\times 3$-matrix we get when we let $\partial$ act on the coefficients of $M$.
By Lemma \ref{derivation} it follows the map
\[ D_{\partial}+\partial(M):A^3\rightarrow A^3 \]
induce a first order differential operator 
\[ \nabla(\partial):\Omega \rightarrow \Omega .\]

Let $\nabla(\partial_1)$ be the following operator:
\[ \nabla(\p)=D_{\p}+\p(M) .\]
One checks that the following holds:
\begin{align}
&\label{PARTIALone} \p(M)=
\end{align}
\[
\begin{pmatrix} -pqx^{p-1}y^{q-1}  & -p(p-1)x^{p-2}y^q +\frac{p^2}{q}x^{2(p-1)}  & -\frac{qp(p-1)}{r}x^{p-2}y^{q-1}z   \\
                -\frac{q^2}{p}y^{2(q-1)}+q(q-1)x^py^{q-2} &  pqx^{p-1}y^{q-1} & \frac{pq(q-1)}{r}x^{p-1}y^{q-2}z \\
 -\frac{qr}{p}y^{q-1}z^{r-1}  & \frac{pr}{q}x^{p-1}z^{r-1}  &  0 
\end{pmatrix}.
\]

One checks that  $\pp(M)$ is the following matrix:
\begin{align}
&\label{PARTIALtwo} \pp(M)=
\end{align}

\[
\begin{pmatrix} -prx^{p-1}z^{r-1} & -\frac{rp(p-1)}{q}x^{p-2}yz^{r-1} & -p(p-1)x^{p-2}z^r+\frac{p^2}{r}x^{2(p-1)} \\
        -\frac{qr}{p}y^{q-1}z^{r-1} & 0 & \frac{pq}{r}x^{p-1}y^{q-1} \\
     -\frac{r^2}{p}z^{2(r-1)}+r(r-1)x^pz^{r-2} & \frac{pr(r-1)}{q}x^{p-1}yz^{r-2} & prx^{p-1}z^{r-1}
\end{pmatrix}.
\]

One checks that $\ppp(M)$ is the following matrix:
\begin{align}
&\label{PARTIALTHREE} \ppp(M)=
\end{align}
\[
\begin{pmatrix} 0 &  -\frac{pr}{q}x^{p-1}z^{r-1} & \frac{pq}{r}x^{p-1}y^{q-1} \\
       -\frac{rq(q-1)}{p}xy^{q-2}z^{r-1} & -qry^{q-1}z^{r-1} & -q(q-1)y^{q-2}z^r+\frac{q^2}{r}y^{2(q-1)} \\
    \frac{qr(r-1)}{p}xy^{q-1}z^{r-2} & -\frac{r^2}{q}z^{2(r-1)}+r(r-1)y^qz^{r-2} & qry^{q-1}z^{r-1}
\end{pmatrix}.
\]

We get well defined first order differential operators
\[ D_{\partial_i}+\partial_i(M):A^3 \rightarrow A^3.\]
for $i=1,2,3$.

\begin{lemma} \label{formone} The following formula holds:
\begin{align}
&\label{done}   (D_{\p}+\p(M))(dF)=(p-q)x^{p-1}y^{q-1}dF. \\
&\label{dtwo}   (D_{\pp}+\pp(M))(dF)=(p-r)x^{p-1}z^{r-1}dF.\\
&\label{dthree} (D_{\ppp}+\ppp(M))(dF)=(q-r)y^{q-1}z^{r-1}dF.
\end{align}
\end{lemma}
\begin{proof} The proof is a straight forward calculation.
\end{proof}

\begin{lemma} \label{DIFFone} The map $\nabla(\partial_i):\Omega \rightarrow \Omega$ defined by
\[ \nabla(\partial_i)(w)=D_{\partial_i}(w)+\partial_i(M)(w) \]
is a well defined differential operator of order one for $i=1,2,3$.
\end{lemma}
\begin{proof} Let $w$ in $\Omega$ be the equivalence class of the element $a=a_1u_1+a_2u_2+a_3u_3\in A^3$.
By definition
\[ \nabla(\partial_i)(w)=\overline{(D_{\partial_i}+\partial_i(M))(a)}.\]
Let $udF\in A^3$. It follows from Lemma \ref{formone} that 
\[ \nabla(\partial_i)(\overline{udF})=u\nabla(\partial_i)(\overline{dF})+\p(u)\overline{dF}=\]
\[  v \overline{dF}=0\]
where $v$ is one of the terms given in the formulas \ref{done}, \ref{dtwo} and \ref{dthree}.
It follows $\nabla(\p)\in \Diff^1(\Omega)$ is a well defined differential operator.
\end{proof}

\begin{theorem} \label{main} The operators $\nabla(\partial_i)$ for $i=1,2,3$ define an algebraic connection
\[ \nabla: \Der_\Z(A)\rightarrow \End_\Z(\Omega) \]
for all integers $p,q,r\geq 2$.
\end{theorem}
\begin{proof} The claim follows from Lemma \ref{DIFFone}.
\end{proof}

Note: By the formulas \ref{PARTIALone}, \ref{PARTIALtwo} and \ref{PARTIALTHREE} it follows the connection $\nabla$ from Theorem \ref{main} 
is defined over any field of characteristic different from $p,q$ and $r$.

Since for any $i,j$ we have by the previous section that the endomorphism  $[\partial_i(M),\partial_j(M)]$
induce the curvature
\[ R_\nabla(\partial_i \wedge \partial_j)\in \End_A(\Omega) \]
it follows the connection $\nabla$ is non flat for any integers $p,q,r\geq 2$: One checks that the $A$-linear endomorphism
\[ [\partial_i(M),\partial_j(M)]: A^3 \rightarrow A^3 \]
induce a non zero endomorphism of $\Omega$.

\begin{corollary} The connection $\nabla$ is non flat in general.
\end{corollary}
\begin{proof} The proof follows from the discussion above.
\end{proof}

The module $\Omega$ is non trivial in general hence it is not clear how to construct a flat algebraic connection on $\Omega$. 

Note: If $A$ is a regular $k$-algebra over an algebraically closed field it follows $\Diff(A)$ is generated by $\Der_k(A)$ as ring.
Hence to give a flat connection 
\[ \nabla: \Der_k(A) \rightarrow \End_k(E) \]
is equivalent to give a left $\Diff(A)$-module structure
\[ \Diff(A)\times E \rightarrow E \]
on $E$ lifting the left $A$-module structure. Given a finitely generated projective $A$-module $E$ the first Chern class
\[ c_1(E)\in \H^2(\Der_k(A),A) \]
is an obstruction to the existence of a flat connection on $E$: If $c_1(E)\neq 0$ it follows $E$ has no flat connections. 

The theory of D-modules is about the study of left modules $\E$ over the sheaf of differential operators $\Diff_X$ where $X$ is a complex manifold,
a differentiable manifold or a smooth algebraic variety. Such $\Diff_X$-modules are usually constructed using Grothendieck's six 
operations for $\Diff_X$-modules. The theory of $\Diff_X$-modules and its relationship to representation theory, the geometry of 
flag varieties and arithmetic is an active area of research (see \cite{borel}, \cite{brylinski}, 
\cite{deligne}, \cite{katz}, \cite{simpson0}, \cite{simpson1} and \cite{simpson2} ).

As indicated above: If one considers a finitely generated projective $A$-module $E$ it follows from Theorem \ref{mainlie} that the connection
$\nabla$ constructed using a projective basis associated to a splitting $s$ is seldom flat. The Lie product
$ [\delta(M), \eta(M)]$ is seldom zero for all derivations $\delta, \eta \in \Der(A)$.

Recall that $\p,\pp$ and $\ppp$ are generators of $\Der_\Z(A)$ as left $A$-module.

\begin{lemma} The following formulas hold in $\Der_\Z(A)$:
\begin{align}
&\label{comone}   [\p,\pp] =p(p-1)x^{p-2}\ppp \\
&\label{comtwo}   [\p,\ppp]=-q(q-1)y^{q-2}\pp \\
&\label{comthree} [\pp,\ppp]=r(r-1)z^{r-2}\p.
\end{align}
\end{lemma}
\begin{proof} The proof is left to the reader.
\end{proof}

Recall that for any derivation $\partial \in \Der_\Z(A)$ 
\[ A_{\partial}=D_{\partial}+\partial(M):A^3\rightarrow A^3 \]
be the corresponding differential operator of order one on $A^3$.
Let moreover
\[ r(\delta \wedge \eta)=[A_{\delta}, A_{\eta}]-A_{[\delta, \eta]}.\]
It follows $r(\delta \wedge \eta)\in \End_A(A^3)$ and there is a commutative diagram
\[
\diagram  A^3 \rto^{r(\delta \wedge \eta)} \dto^p & A^3 \dto^p \\
          E \rto^{R_\nabla(\delta \wedge \eta)} & E
\enddiagram
\]
of $A$-linear maps.

Because of the equality 
\[ tr(R_\nabla(\delta \wedge \eta))=-tr([\delta(M),\eta(M)]|_K) \]
we seek to calculate the map
\[ [\delta(M),\eta(M)]|_K:K\rightarrow K \]
and its trace.
By definition $K$ is the free left $A$-module of rank one on the element $dF$.
We get
\[ [\delta(M),\eta(M)]|_K=r(\delta \wedge \eta)=\]
\[ [A_{\delta},A_{\eta}]-A_{[\delta, \eta]}=\]
\[ A_{\delta}\circ A_{\eta}-A_{\delta} \circ A_{\eta} -A_{[\delta, \eta]}.\]
In the case of $\p$ and $\pp$ we get
\[ r(\p \wedge \pp)=A_{\p}\circ A_{\pp}-A_{\pp}\circ A_{\p}-A_{[\p,\pp]} =\]
\[ A_{\p}\circ A_{\pp}-A_{\pp}\circ A_{\p}-p(p-1)x^{p-2}A_{\ppp}.\]
We calculate
\[ A_{\p}(A_{\pp}(dF)=(p-r)x^{p-2}y^{q-1}z^{r-1}((p-q)x^p+q(p-1))dF.\]
Moreover we get
\[ A_{\pp}(A_{\p}(dF))=(p-q)x^{p-2}y^{q-1}z^{r-1}((p-r)x^p+(p-1)r)dF\]
and
\[ -A_{[\p,\pp]}(dF)=-(q-r)p(p-1)x^{p-2}y^{q-1}z^{r-1}dF.\]
It follows
\[ r(\p \wedge \pp)(dF)\]
\[ A_{\p}\circ A_{\pp}(dF) -A_{\pp}\circ A_{\p}(dF)-p(p-1)x^{p-2}A_{\ppp}(dF)=0.\]
It follows that 
\[ [\p(M), \pp(M)]|_K=0\]
hence $tr([\p(M),\pp(M)]|_K)=0$.

\begin{lemma} The following formulas hold:
\[ r(\partial_i \wedge \partial_j)(dF)=0 \]
for $i,j=1,2,3$.
\end{lemma}
\begin{proof} The proof is done by a calculation similar to the one above.
\end{proof}

\begin{corollary} \label{trace} The following holds:
\[ tr(R_\nabla)=0 .\]
\end{corollary}
\begin{proof} Since $[\partial_i(M),\partial_j(M)]|_K=r(\partial_i \wedge \partial_j)=0$ the result follows from the formula
\[ tr(R_\nabla(\delta \wedge \eta)=-tr([\delta(M),\eta(M)]|_K.\]
\end{proof}

\begin{corollary} \label{zerochern} $c_1(\Omega)=0$ in $\H^2(\Der_\Z(A),A)$ for all integers $p,q,r\geq 2$.
\end{corollary}
\begin{proof} This follows from Corollary  \ref{trace}.
\end{proof}

\begin{example} Chern classes in deRham cohomology and trivial summands.\end{example}

When $A$ is a commutative unital ring and $E$  a rank $d$ projective $A$-module several authors have asked for an obstruction $o(E)$ of $E$
with the property that $o(E)=0$ if and only if $E \cong E' \oplus A$. Recent work of many authors (see \cite{sridharan}) relate the Chern class
$c_d(E)$ in the Chow group of $X:=\Spec(A)$ to this problem. Corollary \ref{zerochern} proves that there is no relation between the  Chern class in deRham cohomology
$c_i(E)$  and this problem.  If $f:=x^2+y^2+z^2-1$ and $A:=k[x,y,z]/(f)$ with $k$ the real numbers, it is well known that the projective module
$T(X):=\Der_k(A)$ does not split off a trivial summand $T(X) \cong L\oplus A$. Still by Corollary \ref{zerochern} it follows $c_i(T(X))=0$ for all $i \geq 1$.
Hence the Chern classes in deRham cohomology does not detect trivial summands.

From Corollary \ref{trace} it is not clear if the module $\Omega$ has a flat connection. Since the first Chern class $c_1(\Omega)$ 
is equal to zero this characteristic class does not give any information on this existence problem. Since $\H^2(\Der_\C(A),A)$ calculates
singular cohomology of the underlying complex manifold $X_\C$ of $X=\Spec(A)$ it follows $c_1(\Omega)$ is a topological obstruction to the existence
of a flat connection on $\Omega$.  

\begin{example} A line bundle with no flat algebraic connections.\end{example}

The following example is related to an example in \cite{loday}. Let $A_{pqr}=\mathbb{C}[x,y,z]/f$ where
$f=x^{2p}+y^{2q}+z^{2r}-1$.  Let $S_{pqr}=\Spec(A_{pqr})$.
In \cite{loday} the author gives an example
of an idempotent for a line bundle $L$ on $S_{111}$ where the image $Ch(L)$ in $\H^*_{\DR}(A_{111})$ is non-zero.
The idempotent is given as follows for arbitrary integers $p,q$ and $r\geq 1$. Let
\[p=
\begin{pmatrix} x^p & y^q+iz^r \\
              y^p-iz^r & -x^p
\end{pmatrix}.
\]
One checks that $p^2=I$ as an element of $\End_{A_{pqr}}(A_{pqr}^2)$. 
It follows the element $M=\frac{1}{2}(p+1)$ is an idempotent for $A_{pqr}^2$ hence if we let
$N=I-M$ we get two endomorphisms of $A_{pqr}^2$ with $M^2=M, N^2=N$ and $M+N=I$. It follows by the results in the first section there is an 
isomorphism
\[ ker(M)\oplus ker(N)\cong A_{pqr}^2 \]
of $A_{pqr}$-modules. Let $K=ker(N), L=ker(M)$.
We get an exact sequence
\[ 0 \rightarrow K \rightarrow A_{pqr}^2 \rightarrow L \rightarrow 0\]
where $L$ is a locally free $A_{pqr}$-module of rank one. Note: $ker(N)=im(M)$ and $ker(M)=im(N)$. It follows there is an isomorphism
\[ L\cong A_{pqr}/im(M)\]
of $A_{pqr}$-modules.
There is an action of the group $G=\Z/p \oplus \Z/q\oplus \Z/r$ 
\[ \sigma: G\times S_{pqr}\rightarrow S_{pqr} .\]
Let $e_j=e^{\frac{2\pi i}{j}}$ and let $u=e_p, v=e_q$ and $w=e_r$.
Define
\[ \sigma((a,b,c), (x,y,z))=(u^ax, v^by, w^cz).\]
It follows $S_{pqr}/G\cong S_{111}$. The canonical map $\pi:S_{pqr}\rightarrow S_{111}$  defined by
\[ \pi(x,y,z)=(x^p, y^q, z^r) \]
is a principal $G$-bundle.
One checks that for any element $(a,b,c)\in \pi^{-1}(x,y,z)$ there is a bijection
\[ \pi^{-1}(x,y,z)\cong \{ (e_p^ka, e_q^lb, e_r^mc) \} \]
where $k=0,\ldots , p-1, l=0,\ldots ,q-1$ and $m=0\ldots , r-1$. 

Since 
\[ Ch:\K_0(A_{111})\rightarrow \H^*_{\DR}(A_{111}) \]
has the property that $Ch_i(L)\neq 0$ for $i\geq 1$ it follows the (non trivial) line bundle $L$ has no flat connections.
The classes $Ch_i(L)$ are independent of choice of connection.
We get a connection on $L$ using the theory in the first section
\[\nabla: \Der_{\mathbb{C}}(A_{pqr})\rightarrow \End_{\mathbb{C}}(L) .\]
Its curvature
\[ R_{\nabla}: \Der_{\mathbb{C}}(A_{pqr})\wedge \Der_{\mathbb{C}}(A_{pqr})\rightarrow \End_{A_{pqr}}(L) \]
is given by the following formula. Let $\delta,\eta\in \Der_{\mathbb{C}}(A_{pqr})$.
\[ R_\nabla(\delta\wedge \eta)=[\delta(M), \eta(M)]. \]
The module $\Der_{\mathbb{C}}(A_{pqr})$ is generated by the following derivations:
\[ D_1=qy^{2q-1}\partial_x-px^{2p-1}\partial_y, D_2=rz^{2r-1}\partial_x-px^{2p-1}\partial_z, D_3=qy^{2q-1}\partial_z-rz^{2r-1}\partial_y.\]

Assume now $p=q=r=1$.

A calculation shows that
\[ D_1(M)=\frac{1}{2}
\begin{pmatrix} y & -x \\
              -x & -y
\end{pmatrix}
\]
and
\[ D_2(M)=\frac{1}{2}
\begin{pmatrix} z & -ix \\
                ix & -z
\end{pmatrix}
\]
and
\[ D_3(M)=\frac{1}{2}
\begin{pmatrix} 0 & z-iy \\
               z+iy & 0
\end{pmatrix}.
\]
We get an algebraic connection
\[ \nabla:\Der_{\mathbb{C}}(A)\rightarrow \End_{\mathbb{C}}(L)\]
defined by
\[ \nabla(D_i)=D_{D_i}+D_i(M)\]
for $i=1,2,3$.
We may calculate the curvature of $\nabla$.
We get
\[ R_\nabla(D_1\wedge D_2)=[D_1(M),D_2(M)]=\]
\[\frac{1}{4}
\begin{pmatrix} -2ix^2 &  2x(z-iy) \\
                -2x(z+iy) &  2ix^2
\end{pmatrix},
\]
\[ R_\nabla(D_1\wedge D_3)=[D_1(M),D_3(M)]=\]
\[\frac{1}{4}
\begin{pmatrix} -2ixy & 2y(z-iy) \\
                -2y(z+iy) &  2ixy
\end{pmatrix}
\]
and
\[ R_\nabla(D_2\wedge D_3)=[D_2(M),D_3(M)]=\]
\[\frac{1}{4}
\begin{pmatrix} -2ixz & 2z(z-iz) \\
               -2z(z+iy) &  2ixz
\end{pmatrix}.
\]

\begin{proposition} \label{nonzerotrace} The following formulas hold:
\begin{align}
\label{tr1}tr(R_\nabla(D_1\wedge D_2))=&-ix \\
\label{tr2}tr(R_\nabla(D_1\wedge D_3))=&-iy \\
\label{tr3}tr(R_\nabla(D_2\wedge D_3))=&-iz 
 \end{align}
\end{proposition}
\begin{proof}  There is a canonical isomorphism
\[ \rho: L^*\otimes_A L\cong \End_A(L) \]
defined by
\[ \rho(x\otimes e)(u)=x(u)e.\]
By definition it follows
\[ L\cong A\{u_1,u_2\}/im(M).\]
Let $e_1$ and $e_2$ be the classes of $u_i$ in $L$. Let $R=R_\nabla(D_1\wedge D_2)$. It follows
\[ tr(R)=x_1(R(e_1))+x_2(R(e_2)).\]
By the calculation above we get
\[ R(e_1)=\frac{1}{4}
\begin{pmatrix}-2ix^2 \\
         -2x(z+iy)
\end{pmatrix}
\]
and
\[ R(e_2)=\frac{1}{4}
\begin{pmatrix} 2x(z-iy)  \\
         2ix^2
\end{pmatrix}.
\]
It follows
\[ tr(R)=\]
\[ -\frac{1}{2}ix^2x_1(e_1)-\frac{1}{2}(xz+ixy)x_1(e_2) +\]
\[ \frac{1}{2}(xz-ixy)x_2(e_1)+\frac{1}{2}ix^2x_2(e_2)=\]
\[ -\frac{1}{2}ix^2(1+x)-\frac{1}{2}(xz+ixy)(y+iz) +\]
\[ \frac{1}{2}(xz-ixy)(y-iz)+\frac{1}{2}ix^2(1-x)=-ix.\]
A similar calculation shows the validity of fomulas \ref{tr2} and \ref{tr2} and the proposition follows.
\end{proof}

\begin{corollary} \label{noflatconn} The projective module $L$ has no flat algebraic connections.
\end{corollary}
\begin{proof} By Proposition \ref{nonzerotrace} it follows $tr(R_\nabla)\neq 0$. By the results in \cite{loday} it follows the Chern class
$c_1(L)$ in $\H^2_{\DR}(A)$ is non zero. Since the class $c_1(L)$ is independent with respect to choice of connection the claim follows.
\end{proof}

\section{On existence of flat algebraic connections}

In this section we relate the results obtained in the previous section to a conjecture on the existence of flat 
connections on complex vector bundles on complex projective manifolds.

\begin{example}\label{noflat}  Projective modules  with no flat connection.\end{example}

Given a connection 
\[ \nabla: \Der_\Z(A) \rightarrow \End_\Z(\Omega) \]
and an $A$-linear map
\[ \phi: \Der_\Z(A) \rightarrow \End_A(\Omega) \]
we get a new connection by letting $\nabla'=\nabla+\phi$. Hence the set of connections on $\Omega$ is parametrized by the set
$\Hom_A(\Der_\Z(A),\End_A(\Omega))$: All other connections on $\Omega$ arise in this way. 
It is not clear if the module $\Omega$ has a flat algebraic connection. One has to investigate if there exists a potential $\phi$
with the property that the connection $\nabla'=\nabla + \phi$ is flat. Let 
\[ r_\phi(\delta \wedge \eta)=[\phi(\delta), \phi(\eta)]-\phi([\delta, \eta]).\]

\begin{lemma} The curvature of the connection $\nabla'=\nabla+\phi$ is as follows:
\[ R_{\nabla'}(\delta \wedge \eta)=R_\nabla(\delta \wedge \eta)+r_\phi(\delta \wedge \eta)+[\nabla(\delta),\phi(\eta)]-[\nabla(\eta),\phi(\delta)].\]
\end{lemma}
\begin{proof} The proof is left to the reader.
\end{proof}

Hence it is not immediate that there exists a potential $\phi:\Der_\Z(A)\rightarrow \End_A(\Omega)$ such that
$R_{\nabla'}=0$. There are many equations that should be fullfilled.

Dually one may consider the connection
\[ \nabla:\Omega \rightarrow \Omega \otimes_A \Omega \]
induced by the projective basis $e_i, x_j$.
The trace of its curvature $tr(R_\nabla)$ gives rise to the first Chern class $c_1(\Omega)\in \H^2_{\DR}(A)$. This class is independent 
with respect to choice of connection. If this class is non-zero it follows $\Omega$ has no flat algebraic connections.

\begin{example} The Chern character of the tangent and cotangent bundle on a smooth hypersurface.\end{example}

Let in the following $f(x_1,\ldots, x_n)\in k[x_1,\ldots, x_n]$ be a non constant polynomial with $k$ the field of real or complex numbers and let $A=k[x_1,\ldots, x_n]/f$. Assume
$X=\Spec(A)$ is  regular. Let
\[ df=f_{x_1}dx_1+\cdots f_{x_n}dx_n\in A\{dx_1,\ldots ,dx_n\}.\]
It follows $df=0$ if and only if $f$ is a constant. We get a split exact sequence
\[ 0\rightarrow Adf \rightarrow A\{dx_1,\ldots dx_n\} \rightarrow \Omega:=\Omega^1_A\rightarrow 0\]
of left $A$-modules. Choose a splitting $s$ and a projective basis $\{B,B^*\}$ for the above sequence, with corresponding connections 

 \[ \nabla'_B: \Der_k(A) \rightarrow \End_k(\Omega) \]

and

\[ \nabla_B: \Omega^1 \rightarrow \Omega^1 \otimes \Omega^1.\]

It follows we may define the $i$'th Chern character

\[  ch_i(\Omega):= \frac{1}{i!}tr(R_{{\nabla}_B}^i) \in \H^{2i}_{DR}(A)\]

where $\H^j_{DR}(A)$ is algebraic deRham cohomology of $A$. Let $T_A:=\Der_k(A)$.

\begin{lemma} \label{zero} The following holds: $ch_i(\Omega^1_A)=0$ and $ch_i(T_A)=0$ in $\H^{2i}_{DR}(A)$ for all $i \geq 1$.
\end{lemma}
\begin{proof}  There is a well defined map of abelian groups - the Chern character

\[ Ch: \K_0(A) \rightarrow \H^{2*}_{DR}(A) \]

and in $\K_0(A)$ there is an equality

\[  [\Omega^1_A]=(n-1)[A].\]

The Chern character is independent with respect to choice of connection and since the free rank one $A$-module $A$
has a flat connection the claim follows. A similar argument proves the same formula for $ch_i(T_A)$. The Lemma follows from this.
\end{proof}

Note: In Lemma \ref{zero} $\H^*_{DR}(A)$ is defined using the module of Kahler differentials $\Omega^1_A$, its exterior powers
$\wedge^p \Omega^1_A$ and the "naive" algebraic deRham complex. In general when $A$ is non-regular one uses a compactification 
to define algebraic deRham cohomology (see \cite{hartshorne}).

\begin{example} Regular hypersurfaces, affine Calabi-Yau manifolds and Suslins theorem on stably trivial projective modules. \end{example}

Suslin proved in the 1970's (see \cite{suslin}) that if the base field $k$ is algebraically closed and if $krdim(A)=n$, any stably free rank $n$ projective $A$-module is free. Hence if $A$ is a regular ring of finite type over an algebraically closed field $k$, it follows the tangent and cotangent bundles are free. Lemma \ref{zero} holds over any field $k$, and over an arbitrary field $k$, the tangent and cotangent bundle of a smooth affine hypersurface $X:=V(f)$ will be non-trivial. Hence Lemma \ref{zero} proves that when $A$ is regular and of finite type over a field, the Chern class in "naive" algebraic deRham cohomology and Lie-Rinehart cohomology is a "weak" invariant. It does not detect that the tangent and cotangent bundle is non-trivial in this case. 

Note moreover that the above construction proves that if $k$ is the field of complex numbers  it follows the cotangent and tangent bundle $\Omega^1_X$ and $T_X$ are trivial. This follows from the Suslin theorem since $\Omega^1_X$ is stably trivial and since $T_X$ is the dual of the cotangent bundle. 

\begin{theorem} Let $k$ be the field of complex numbers and let $X:=V(f) \subseteq \mathbb{A}^n_k$ be a regular hypersurface. It follows the underlying complex manifold of $X$
is a Calabi-Yau manifold.
\end{theorem}
\begin{proof} By the above argument it follows the canonical bundle $\omega_X$ is trivial, hence $X$ is a Calabi-Yau manifold. 
\end{proof}

Note: Usually a Calabi-Yau manifold is defined to be compact Kahler manifold with trivial canonical bundle and an affine hypersurface is not compact in general.


By Theorem \ref{main} and Lemma \ref{trace} it follows the cotangent bundle $\Omega^1_A$  on $f:=x^p+y^q+z^r-1$ has a non flat algebraic connection
\[ \nabla: \Der_{k}(A)\rightarrow \End_A(\Omega^1_A).\]

Here $A:=k[x,y,z]/(f)$. It is not clear how to construct a potential $\phi$ such that $\nabla'=\nabla+\phi$ is a flat connection on $\Omega^1_A$, hence it is not clear if $\Omega^1_A$ has a flat algebraic connection. 

\begin{example} Non-abelian extensions on the real 2-sphere and torsors . \end{example}

If the field $k$ is arbitrary, the cotangent bundle $\Omega^1_A$ will be non-trivial in general - if $k$ is the real numbers and $p=q=r=2$ we get the dual of the tangent bundle of the real 2-sphere $S(2):=\Spec(A)$ and this is a topologically non-trivial vector bundle. Let $T_{S(2)}$ be the tangent bundle of $S(2)$ and let $\Omega^1_{S(2)}$ be the cotangent bundle. 
By the calculation in Theorem \ref{main} we get a non-flat algebraic connection

\[ \nabla: T_{S(2)} \rightarrow \End_{\mathcal{O}_{S(2)}}( \Omega^1_{S(2)}) ,\]
with $tr(K_{\nabla})=0$. This and Lemma \ref{zero} implies $c_i(\Omega^1_{S(2)})=0$ for all  $i\geq 1$.

There is a "refined "characteristic class

\[ c(\Omega^1_A) \in \Ext^1(\Der_k(A), \End_A(\Omega^1_A)) ,\]

and the class $c(\Omega^1_A)$ is defined using the connection $\nabla$. The class is independent of choice of connection 
and $c(\Omega^1_A)$ is trivial if and only if $\Omega^1_A$ has a flat algebraic connection (see \cite{maa145}, Theorem 2.15).
The set  $ \Ext^1(\Der_k(A), \End_A(\Omega^1_A))$ is the set of non-abelian extensions of $\Der_k(A)$ with $\End_A(\Omega^1_A)$. It is a pointed set (see \cite{maa145}).  By Lemma \ref{noflatconn} it follows the class $c(-)$ is a non-trivial characteristic class. 

Note: The pointed cohomology set $\Ext^1(\Der_k(A), \End_A(\Omega^1_A))$ is a \emph{torsor} on a certain cohomology group defined using the Lie-Rinehart cohomology of a connection. 
Hence the pointed cohomology set has additional structure.

\begin{example} \label{existence} Complex projective ellipsoid surfaces. \end{example}

Consider the polynomial
\[ F= X^pT^{q+r}+Y^qT^{p+r}+Z^rT^{p+q}-T^{p+q+r}\in \C[X,Y,Z,T].\]
Let $B=\C[X,Y,Z,T]/F$.
We get a complex projective surface
\[ S=\Proj(B)\subseteq \P^3_{\C}.\]
Let $U=D(T)\cap S$. It follows $U$ is the affine ellipsoid surface $\Spec(A)$ studied above. 

\begin{lemma} \label{counter} Assume $\Omega$ has no flat algebraic connections and assume the Atiyah 
class $a(\Omega^1_S)=0$. It follows $\Omega^1_S$ is a complex vector bundle on $S$ with a connection.
It moreover follows $\Omega^1_S$ does not have a flat connection.
\end{lemma}
\begin{proof} Assume $\Omega^1_S$ has a flat connection
\[ \nabla: \Der_{\C}(\O_S)\rightarrow \End_{\C}(\Omega^1_S).\]
We may restrict $\nabla$ to $U=\Spec(A)$ to get a flat connection
\[ \nabla|_U:\Der_{\C}(A)\rightarrow \End_{\C}(\Omega^1_S|_U) .\]
By functoriality $\Omega^1_S|_U=\Omega$.
By the above results it follows $\Omega$ does not have a flat connection which is a contradiction.
\end{proof}

Hence if $\Omega$ has no flat algebraic connections and if the Atiyah class $a(\Omega^1_S)$ is zero we get by 
Lemma \ref{counter} a counterexample to an old conjecture 
(see \cite{atiyah}) on existence of flat connections for finite rank vector bundles on complex projective manifolds.
The conjecture says that if a complex vector bundle on a complex projective manifold has a holomorphic  connection then it has a flat 
holomorphic connection. Note that for a complex vector bundle $\E$ on a complex projective manifold $X$ a holomorphic connection
is given as a left splitting of the canonical map
\[ p: \J^1(\E)\rightarrow \E \]
constructed in \cite{maa2}. Since $X$ is projective it follows from \cite{serre} $X$ is a smooth algebraic variety.
When $\E$ is locally free it follows $\J^1(\E)$ is locally free hence a holomorphic connection $\nabla$ is algebraic since it is a left
splitting of $p$ which is a surjective map of algebraic vector bundles.

\begin{example} Locally free finite rank sheaves on complex projective manifolds.\end{example}

Assume $X\subseteq \mathbb{P}^n_\C$ is a complex projective manifold. It follows 
that $X$ is a smooth algebraic  variety. Let $\E$ be a finite rank locally free sheaf on $X$ and assume
\[ \nabla: \E\rightarrow \Omega^1_X\otimes \E \]
is a holomorphic connection. By the above arguments it follows $\nabla$ is algebraic hence it is locally defined by polynomial functions.
Let $U=\Spec(A)\subseteq X$ be an affine open subset and let $E=\E|_U$ be the induced $A$-module. We get
an induced connection
\[ \nabla|_U: E\rightarrow \Omega^1_A \otimes E  .\]
If the corresponding Chern class $c_1(E)$ in $\H^2(\Der_\Z(A),A)$ or $\H^2_{\DR}(A)$ is non zero we get a counter example to the conjecture.
The class $c_1(E)$ is independent with respect to choice of connection  hence $E$ has no flat connections. By the calculation in the 
previous section it seems quite unlikely that $c_1(E)=0$ in $\H^2(\Der_\Z(A),A)$ for every open affine subset $\Spec(A)\subseteq X$.

By Lemma \ref{zero} this argument
does not work for smooth hypersurfaces in $\mathbb{P}^n_\C$.  

The class $c_1(E)$ may be calculated using Groebner bases.

If there is a map of $\C$-vector spaces
\[ f:\Ext^1_A(E,\Omega^1_A \otimes E)\rightarrow \H^2_{\DR}(A) \]
with 
\[ f(a(E))=c_1(E) \]
this argument fails.

\end{document}